# On the asymptotic joint distribution of sample space–time covariance estimators

BO LI[1], MARC G. GENTON[2] and MICHAEL SHERMAN[3]

[1]*NCAR, 1850 Table Mesa Dr., Boulder, CO 80305, USA. E-mail: boli@ucar.edu*
[2]*Department of Econometrics, University of Geneva, Bd du Pont-d'Arve 40, CH-1211 Geneva 4, Switzerland. E-mail: Marc.Genton@metri.unige.ch and Department of Statistics, Texas A&M University, College Station, TX 77843-3143, USA. E-mail: genton@stat.tamu.edu*
[3]*Department of Statistics, Texas A&M University, College Station, TX 77843-3143, USA. E-mail: sherman@stat.tamu.edu*

We study the asymptotic joint distribution of sample space–time covariance estimators of strictly stationary random fields. We do this without any marginal or joint distributional assumptions other than mild moment and mixing conditions. We consider several situations depending on whether the observations are regularly or irregularly spaced and whether one part or the whole domain of interest is fixed or increasing. A simulation experiment illustrates the theoretical results.

## 1. Introduction

Let $\{Z(\mathbf{x}) : \mathbf{x} \in D \subset \mathbb{R}^d\}$, $d \geq 1$, be a strictly stationary random field with covariance function

$$C(\mathbf{k}) = \text{cov}\{Z(\mathbf{x}), Z(\mathbf{x} + \mathbf{k})\},$$

where $\mathbf{k}$ denotes an arbitrary lag in $\mathbb{R}^d$. Let $\mathbf{\Lambda}$ be a set of lags and let $m$ denote the cardinality of $\mathbf{\Lambda}$. Define $\mathbf{G} = \{C(\mathbf{k}) : \mathbf{k} \in \mathbf{\Lambda}\}$ to be the length $m$ vector of covariances at lags in $\mathbf{\Lambda}$. Let $\widehat{C}_n(\mathbf{k})$ denote a sample-based moment or kernel estimator of $C(\mathbf{k})$ based on observations in a sequence of increasing index sets $D_n \subset D$ and let $\widehat{\mathbf{G}}_n = \{\widehat{C}_n(\mathbf{k}) : \mathbf{k} \in \mathbf{\Lambda}\}$ denote the estimator of $\mathbf{G}$. We are interested in the asymptotic distribution of $\widehat{\mathbf{G}}_n$. This distribution is important for several reasons, such as assessing the covariance structure of a random field. For example, Li, Genton and Sherman (2007) proposed a testing approach for various properties of space–time covariance functions based on the asymptotic distribution of covariance estimators in a particular setting. For the particular case of testing directional properties of a random field, Guan, Sherman and Calvin (2004) and Lu and Zimmerman (2001) derived the asymptotic distribution of the spatial variogram when $d = 2$. Surprisingly, the distribution of sample covariances has been investigated only in particular situations in the literature. Under the assumption of a Poisson process for modeling the observations' locations, Masry (1983) proved the asymptotic joint normality of sample autocovariances for time series and Karr (1986) generalized this result to a random field. Brockwell and Davis (1991), page 229, and Fuller (1996), page 333, derived the asymptotic joint normality of sample autocovariances under mild assumptions for stationary time series. For spatially indexed observations, several authors (e.g., Bolthausen (1982), Guyon (1995)) have proven the asymptotic normality of the scalar sample mean under a variety of mixing and moment conditions.

Recognizing various situations occurring in practice, we consider several data structures and asymptotic regimes, depending on whether the observations are regularly spaced or irregularly spaced and whether a subset of dimensions is fixed or the whole domain of interest is increasing. We also allow one dimension to denote time. We make no assumptions on the marginal or the joint distribution of observations other than mild moment and mixing conditions on the random field.

To formally state the asymptotic properties of the vector of sample covariances $\widehat{\mathbf{G}}_n$, we need to quantify the strength of dependence in the random field, taking account of different types of spacing of observations. Following Rosenblatt (1956), we define the strong mixing coefficients for a random field with regularly spaced observations as

$$\alpha_b(r) = \sup_{E_1, E_2} \{|P(A_1 \cap A_2) - P(A_1)P(A_2)| : A_i \in \mathfrak{F}(E_i), |E_i| \le b, \\ i = 1, 2, d(E_1, E_2) \ge r\}, \tag{1.1}$$

where $|E|$ denotes the cardinality of the set $E$, $\mathfrak{F}(E)$ denotes the $\sigma$-algebra generated by the random variables $\{Z(\mathbf{x}) : \mathbf{x} \in E\}$ and where $d(E_1, E_2) = \inf\{\sup_j |\mathbf{x}_{1j} - \mathbf{x}_{2j}| : \mathbf{x}_1 \in E_1, \mathbf{x}_2 \in E_2, j = 1, \ldots, d\}$. The supremum in (1.1) is taken over all compact and convex subsets $E_1 \subset \mathbb{R}^d$ and $E_2 \subset \mathbb{R}^d$ such that $d(E_1, E_2) \ge r$. For a random field with irregularly spaced observations, we define the mixing coefficients following Politis, Paparoditis and Romano (1998). This definition, denoted (1.1′), is formed by imposing an additional condition that $E_2$ is a shift of $E_1$ in (1.1). Note that $|E|$ denotes the Lebesgue measure (volume) of $E$ in (1.1′).

If the observations are independent, then $\alpha_b(r) = 0$ for all $r > 0$. We need $\alpha_b(r)$ to approach 0 for large $r$, at some rate depending on $b$. In order to appropriately define this rate, we decompose $D_n$ into $D_n = \mathcal{F} \times \mathcal{I}_n$, where $\mathcal{F} \subset \mathbb{R}^p$, $\mathcal{I}_n \subset \mathbb{R}^q$ and $p + q = d$. Suppose that $\mathcal{F}$ is a fixed space, in the sense that finitely many observations are located within this space, and $\mathcal{I}_n$ is an increasing space. Following Sherman and Carlstein (1994), we assume the mixing condition for $\mathcal{I}_n$

$$\sup_b \frac{\alpha_b(r)}{b} = \mathrm{O}(r^{-\varepsilon}) \qquad \text{for some } \varepsilon > q. \tag{C1}$$

From Proposition 3.1 of Künsch (1982), it can be shown that (C1) holds for a class of Gibbs fields. Bradley (1993) has also studied and justified the use of mixing condition

(C1). We account for the shape of the domain in which we observe data as in Sherman (1996). Let $\mathcal{A}$ denote the interior of a closed hypersurface contained in a $q$-dimensional hypercube with edge length 1. $\mathcal{A}_n$ is the inflation of $\mathcal{A}$ by a factor $n$. We define $\mathcal{I}_n = \mathbb{Z}^q \cap \mathcal{A}_n$ if the observations are regularly spaced, and $\mathcal{I}_n = \mathcal{A}_n$ otherwise. Note that $\mathcal{I}_n$ satisfies

$$|\mathcal{I}_n| = \mathrm{O}(n^q), \qquad |\partial \mathcal{I}_n| = \mathrm{O}(n^{q-1}), \tag{C2}$$

where $\partial \mathcal{I}_n$ denotes the boundary of $\mathcal{I}_n$. This allows for a wide variety of domains.

If a random field provides sufficient information so that we can estimate its covariance at any arbitrary lag, this random field is said to exhibit *continuous lags*. We use a kernel estimator to estimate covariances in this case. For example, if the observations are irregularly spaced in the subspace $\mathbb{R}^{q_1}$ of $\mathbb{R}^q$ ($q_1 \leq q$), covariances for any lag value in $\mathbb{R}^{q_1}$ can be estimated by a kernel estimator; see Section 3. However, if covariances at only discrete lags are estimable, this random field is said to exhibit *discrete lags* and we use a moment estimator to estimate the covariances in this situation.

Throughout the paper, we assume that the mean of $Z$ is known and equal to 0. If we remove this assumption, let $\widehat{C}_n^*(\mathbf{k})$ and $\widehat{\mathbf{G}}_n^*$ denote the mean corrected estimators of $C(\mathbf{k})$ and $\mathbf{G}$, respectively. In Lemma A.6, we show that $\widehat{\mathbf{G}}_n^*$ and $\widehat{\mathbf{G}}_n$ have the same asymptotic properties. Let $\widehat{C}_n(\mathbf{k})$ denote a moment estimator or a kernel estimator of the covariance $C(\mathbf{k})$ under the zero-mean assumption and let $\lambda_n$ be the bandwidth of the kernel estimator defined over an irregularly spaced subspace $\mathbb{R}^{q_1}$ of $\mathbb{R}^q$. We assume the following moment condition for the moment estimator:

$$\sup_n \mathrm{E}\{|\sqrt{|\mathcal{I}_n|}\{\widehat{C}_n(\mathbf{k}) - C(\mathbf{k})\}|^{2+\delta}\} \leq C_\delta \qquad \text{for some } \delta > 0, C_\delta < \infty. \tag{C3}$$

The moment condition (C3$'$) for the kernel estimator can be obtained from (C3) by simply replacing $|\mathcal{I}_n|$ with $|\mathcal{I}_n|\lambda_n^{q_1}$ and $C(\mathbf{k})$ with $\mathrm{E}\{\widehat{C}_n(\mathbf{k})\}$. The moment condition is only slightly stronger than the existence of the (standardized) asymptotic variance of $\widehat{C}_n(\mathbf{k})$.

In this article, we derive the asymptotic joint distribution of sample covariances for space–time random fields $\{Z(\mathbf{s}, t) : \mathbf{s} \in \mathbb{R}^2, t \in \mathbb{R}\}$. However, the results for this $\mathbb{R}^2 \times \mathbb{R}$ space can be easily extended to the $\mathbb{R}^d$ ($d > 3$) space. Note that we consider an increasing-domain asymptotic framework. For a recent comparison of infill and increasing-domain asymptotics and the consequences for maximum likelihood estimators of covariance parameters, see Zhang and Zimmerman (2005).

The rest of the paper is organized as follows. Section 2 demonstrates the asymptotic joint normality of sample space–time covariances for discrete lags, while Section 3 addresses the distributional behavior in the case of continuous lags. These two sections consider space–time data structures which are common in applications. Section 4 presents a simulation experiment. Appendix A states and proves some useful lemmas. Appendix B contains proofs of all theorems.

## 2. Discrete lags

### 2.1. Regularly spaced observations with an increasing spatio-temporal domain

Consider a strictly stationary mean zero spatio-temporal random field $\{Z(\mathbf{s},t):\mathbf{s}\in\mathbb{R}^2, t\in\mathbb{R}\}$. Let $D_n = S_n \times T_n$ be a finite set of lattice points in $\mathbb{Z}^2 \times \mathbb{Z}$ at which observations are taken. We allow the lattice in $S_n$ and $T_n$ to be defined in different metrics. Let $\mathbf{h}$ denote a lag in space and $u$ denote a lag in time. The classical estimator of the covariance, that is, the sample covariance, is given by

$$\widehat{C}_n(\mathbf{h},u) = \frac{1}{|D_n(\mathbf{h},u)|} \sum_{D_n(\mathbf{h},u)} Z(\mathbf{s},t)Z(\mathbf{s}+\mathbf{h},t+u),$$

where the sum is over $D_n(\mathbf{h},u) = \{(\mathbf{s},t) : (\mathbf{s},t) \in D_n, (\mathbf{s}+\mathbf{h}, t+u) \in D_n\}$ and $|D_n(\mathbf{h},u)|$ denotes the number of distinct elements in $D_n(\mathbf{h},u)$.

We define the strong mixing coefficients of this random field as in (1.1) and assume the mixing condition, boundary condition and moment condition as in (C1), (C2) and (C3) setting $p=0$, $q=3$ and $\mathcal{I}_n = D_n$.

**Theorem 1.** *Let $\{Z(\mathbf{s},t) : \mathbf{s} \in \mathbb{R}^2, t \in \mathbb{R}\}$ be a strictly stationary mean zero spatio-temporal random field observed at lattice points in $D_n \subset \mathbb{Z}^3$ satisfying condition* (C2). *Assume that*

$$\sum_{\mathbf{s}\in\mathbb{Z}^2}\sum_{t\in\mathbb{Z}} |\operatorname{cov}\{Z(\mathbf{0},0)Z(\mathbf{h}_1,u_1), Z(\mathbf{s},t)Z(\mathbf{s}+\mathbf{h}_2,t+u_2)\}| < \infty \qquad (2.1)$$

*for all finite $\mathbf{h}_1, \mathbf{h}_2, u_1$ and $u_2$. Then, $\Sigma = \lim_{n\to\infty}|D_n|\operatorname{cov}(\widehat{\mathbf{G}}_n, \widehat{\mathbf{G}}_n)$ exists, the $(i,j)$th element of which is*

$$\sum_{\mathbf{s}\in\mathbb{Z}^2}\sum_{t\in\mathbb{Z}} \operatorname{cov}\{Z(\mathbf{0},0)Z(\mathbf{h}_i,u_i), Z(\mathbf{s},t)Z(\mathbf{s}+\mathbf{h}_j,t+u_j)\}.$$

*If we further assume that $\Sigma$ is positive definite and that conditions* (C1) *and* (C3) *hold, then $\sqrt{|D_n|}(\widehat{\mathbf{G}}_n - \mathbf{G}) \xrightarrow{d} N_m(\mathbf{0}, \Sigma)$ as $n \to \infty$.*

### 2.2. Observations with a fixed spatial domain and an increasing temporal domain

In many situations, the observations are taken from a fixed space $S \subset \mathbb{R}^2$ at regularly spaced times $T_n$. Let $S(\mathbf{h}) = \{\mathbf{s} : \mathbf{s} \in S, \mathbf{s} + \mathbf{h} \in S\}$ and $|S(\mathbf{h})|$ be the number of elements in $S(\mathbf{h})$. In this situation, we have $p=2$, $q=1$, $\mathcal{F} = S$ and $\mathcal{I}_n = T_n$ in the notation of the Introduction.

In this particular case, we define the mixing coefficient (e.g., Ibragimov and Linnik (1971), page 306) only in the time direction since the domain is increasing only in this direction:

$$\alpha(u) = \sup_{A,B}\{|P(A \cap B) - P(A)P(B)|, A \in \mathfrak{F}_{-\infty}^0, B \in \mathfrak{F}_u^\infty\},$$

where $\mathfrak{F}_{-\infty}^0$ is the $\sigma$-algebra generated by the past-time process until $t = 0$ and $\mathfrak{F}_u^\infty$ is the $\sigma$-algebra generated by the future-time process from $t = u$. We assume the mixing coefficient $\alpha(u)$ satisfies the strong mixing condition:

$$\alpha(u) = \mathrm{O}(u^{-\varepsilon}) \qquad \text{for some } \varepsilon > 0, \tag{2.2}$$

and we also observe that

$$T_n = \{1, \ldots, n\}, \qquad |\partial T_n| = 2 = \mathrm{O}(1). \tag{2.3}$$

Condition (2.3) is a special case of (C2), while condition (2.2) holds for a variety of temporal processes. For example, (2.2) holds for AR(1) processes with normal, double exponential or Cauchy errors by the results of Gastwirth and Rubin (1975). We further assume the moment condition as in (C3) and define the estimator of $C$ as

$$\widehat{C}_n(\mathbf{h}, u) = \frac{1}{|S(\mathbf{h})||T_n|} \sum_{S(\mathbf{h})} \sum_{t=1}^{n-u} Z(\mathbf{s},t) Z(\mathbf{s}+\mathbf{h}, t+u).$$

**Theorem 2.** *Let $\{Z(\mathbf{s},t), \mathbf{s} \in \mathbb{R}^2, t \in \mathbb{R}\}$ be a strictly stationary mean zero spatio-temporal random field observed in $D_n = S \times T_n$, where $S \subset \mathbb{R}^2$ and $T_n$ satisfies condition (2.3). Assume that*

$$\sum_{t \in \mathbb{Z}} |\mathrm{cov}\{Z(\mathbf{0},0)Z(\mathbf{h}_1, u_1), Z(\mathbf{s},t)Z(\mathbf{s}+\mathbf{h}_2, t+u_2)\}| < \infty \tag{2.4}$$

*for all $\mathbf{h}_1 \in S, \mathbf{h}_2 \in S, \mathbf{s} \in S(\mathbf{h}_2)$ and all finite $u_1, u_2$. Then, $\Sigma = \lim_{n\to\infty} |T_n| \mathrm{cov}(\widehat{\mathbf{G}}_n, \widehat{\mathbf{G}}_n)$ exists, the $(i,j)$th element of which is*

$$\frac{1}{|S(\mathbf{h}_i)||S(\mathbf{h}_j)|} \sum_{\mathbf{s}_1 \in S(\mathbf{h}_i)} \sum_{\mathbf{s}_2 \in S(\mathbf{h}_j)} \sum_{t \in \mathbb{Z}} \mathrm{cov}\{Z(\mathbf{s}_1, 0)Z(\mathbf{s}_1 + \mathbf{h}_i, u_i), Z(\mathbf{s}_2, t)Z(\mathbf{s}_2 + \mathbf{h}_j, t + u_j)\}.$$

*If we further assume that $\Sigma$ is positive definite and that conditions (2.2) and (C3) hold, then $\sqrt{|T_n|}(\widehat{\mathbf{G}}_n - \mathbf{G}) \xrightarrow{d} \mathrm{N}_m(\mathbf{0}, \Sigma)$ as $n \to \infty$.*

In Theorem 2, we allow the observations to be either regularly spaced or irregularly spaced in $S$. However, even for irregularly spaced observations, we consider only the covariances of observed spatial lags due to the limited number of observations in $S$. Note that in this section, we require the observations to be taken at the same spatial locations over time, which is very common for monitoring stations. For example, the Irish wind data of Haslett and Raftery (1989) recently analyzed by Gneiting (2002), de Luna and Genton (2005) and Stein (2005) consists of time series of daily average wind speed at eleven meteorological stations in Ireland. In particular, Li, Genton and Sherman (2007) developed space–time symmetry and separability tests based on the asymptotic normality of sample covariance estimators under this specific data structure.

## 3. Continuous lags

### 3.1. Spatially irregularly spaced observations with an increasing spatio-temporal domain

If the observations are spatially irregularly spaced in an increasing domain, we can then estimate the covariance for any spatial lag by employing kernel smoothing. Again, consider a strictly stationary mean zero random field $\{Z(\mathbf{s},t): \mathbf{s} \in \mathbb{R}^2, t \in \mathbb{R}\}$. Let $D_n \subset \mathbb{R}^3$ denote the domain of interest in which observations are taken. We decompose $D_n$ into $D_n = S_n \times T_n$, where $S_n$ is an increasing index set, $S_n \subset \mathbb{R}^2$ and $T_n = \{1, \ldots, n\}$. We view the spatial locations at which $Z$ is observed in $S_n$ as *random* in number and location; specifically, they are generated from a homogeneous 2-dimensional Poisson process with intensity parameter $\nu$.

Let $N$ denote the random point process and $N(B)$ denote the random number of points of $N$ contained in $B$, where $B$ is any given Borel set. We further assume $N$ to be independent of $Z$. We construct an estimator of covariance based on kernel smoothing.

Let $|S_n|$ denote the Lebesgue measure (not the cardinality) of $S_n$ and let $w(\cdot)$ be a bounded, symmetric density function on $\mathbb{R}^2$. Here and henceforth, we use $d\mathbf{s}$ to denote an infinitesimally small disc centered at $\mathbf{s}$. Define $N^{(2)}(d\mathbf{s}_1, d\mathbf{s}_2) \equiv N(d\mathbf{s}_1)N(d\mathbf{s}_2)I(\mathbf{s}_1 \neq \mathbf{s}_2)$, where $I(\mathbf{s}_1 \neq \mathbf{s}_2) = 1$ if $\mathbf{s}_1 \neq \mathbf{s}_2$ and 0 otherwise. The kernel covariance estimator over $D_n$ is given by

$$\widehat{C}_n(\mathbf{h}, u) = \frac{1}{\nu^2 |T_n||S_n|} \sum_{t=1}^{n-u} \int_{S_n} \int_{S_n} w_n(\mathbf{h} - \mathbf{s}_1 + \mathbf{s}_2) Z(\mathbf{s}_1, t) Z(\mathbf{s}_2, t+u) N^{(2)}(d\mathbf{s}_1, d\mathbf{s}_2),$$

where $w_n(\mathbf{x}) = \frac{1}{\lambda_n^2} w(\frac{\mathbf{x}}{\lambda_n})$ and $\lambda_n$ is a sequence of positive constants satisfying condition (3.1) below. Here, $\nu$ is assumed to be known, otherwise it can be consistently estimated by $\widehat{\nu} = \frac{N(S_n)}{|S_n|}$.

We define the mixing coefficients as in (1.1′) and assume the mixing and moment conditions as in (C1) and (C3′), setting $p = 0$, $q = 3$, $q_1 = 2$ and $\mathcal{I}_n = D_n$. In addition to the boundary condition in (C2), we need to consider the choice of bandwidth. Specifically, we assume that

$$\lambda_n \to 0, \qquad \lambda_n^2 |S_n| \to \infty. \tag{3.1}$$

Define the following fourth order cumulant function (e.g., Karr (1986)):

$$Q(\mathbf{x}_1, \mathbf{x}_2, \mathbf{x}_3) = \mathrm{E}\{Z(\mathbf{0})Z(\mathbf{x}_1)Z(\mathbf{x}_2)Z(\mathbf{x}_3)\} - C(\mathbf{x}_1)C(\mathbf{x}_3 - \mathbf{x}_2) \\ - C(\mathbf{x}_2)C(\mathbf{x}_3 - \mathbf{x}_1) - C(\mathbf{x}_3)C(\mathbf{x}_2 - \mathbf{x}_1), \quad (3.2)$$

where the $\mathbf{x}$'s denote generic locations and $\mathbf{x} = (\mathbf{s}, t)$ in the spatio-temporal random field.

The following theorem states that $\widehat{C}_n(\mathbf{h}, u)$ is a consistent estimator of $C(\mathbf{h}, u)$ and $\widehat{\mathbf{G}}_n$ is asymptotically jointly normal under some conditions.

**Theorem 3.** *Let* $\{Z(\mathbf{s}, t) : \mathbf{s} \in \mathbb{R}^2,\ t \in \mathbb{R}\}$ *be a strictly stationary mean zero spatio-temporal random field observed in* $D_n = S_n \times T_n$, *where* $S_n \subset \mathbb{R}^2$ *and* $T_n = \{1, \ldots, n\}$. $D_n$ *satisfies condition* (C2). *Assume that the locations* $\mathbf{s}$ *are generated by a homogeneous Poisson process in* $\mathbb{R}^2$ *and that*

$$\sup_u \int_{\mathbb{R}^2} C(\mathbf{h}, u) \, d\mathbf{h} < \infty,$$

$$\sup_{\mathbf{h}, u, u'} \sum_{t \in \mathbb{Z}} \mathrm{E}\{Z(\mathbf{0}, 0)Z(\mathbf{h}, u)Z(\mathbf{0}, t)Z(\mathbf{h}, t + u')\} < \infty,$$

$$\sup_{\mathbf{s}_1, \mathbf{s}_3} \sup_{t_1, t_2, t_3} \int_{\mathbf{s}_2 \in \mathbb{R}^2} |Q\{(\mathbf{s}_1, t_1), (\mathbf{s}_2, t_2), (\mathbf{s}_2 + \mathbf{s}_3, t_2 + t_3)\}| \, d\mathbf{s}_2 < \infty.$$

*Then,* $\mathrm{E}\{\widehat{C}_n(\mathbf{h}, u)\} \to C(\mathbf{h}, u)$ *and* $\Sigma = \lim_{n \to \infty} |T_n||S_n|\lambda_n^2 \operatorname{cov}(\widehat{\mathbf{G}}_n, \widehat{\mathbf{G}}_n)$ *exists, the* $(i, j)$*th element of which is*

$$\frac{1}{\nu^2} \int_{\mathbb{R}^2} w^2(\mathbf{x}) \, d\mathbf{x} \sum_{t \in \mathbb{Z}} \mathrm{E}[Z(\mathbf{0}, 0)Z(\mathbf{h}_i, u_i)\{Z(\mathbf{0}, t)Z(\mathbf{h}_i, t + u_j)I(\mathbf{h}_i = \mathbf{h}_j) \\ + Z(\mathbf{h}_i, t)Z(\mathbf{0}, t + u_j)I(\mathbf{h}_i = -\mathbf{h}_j)\}],$$

*where* $I(\mathbf{u} = \mathbf{v}) = 1$ *if* $\mathbf{u} = \mathbf{v}$ *and 0 otherwise. If we further assume that* $\Sigma$ *is positive definite and the conditions* (C1), (C3′) *and* (3.1) *hold, then*

$$\sqrt{|T_n||S_n|}\lambda_n\{\widehat{\mathbf{G}}_n - \mathrm{E}(\widehat{\mathbf{G}}_n)\} \xrightarrow{d} \mathrm{N}_m(\mathbf{0}, \Sigma)$$

*as* $n \to \infty$.

## 3.2. Irregularly spaced observations with an increasing spatio-temporal domain

In Section 3.1, we discussed the properties of the kernel covariance estimator if the observations are irregularly spaced in $S_n$ and regularly spaced in $T_n$. Another case occurs when the observations are irregularly spaced in the whole space $S_n \times T_n$. Karr (1986) showed mean square consistency of the estimator $\widehat{C}$ of the covariance function $C$ of the random field $Z$ defined on $\mathbb{R}^d$ ($d \geq 2$). It is an extension of the results in Masry (1983) who investigated the covariance estimator for time series, that is, for random fields in $\mathbb{R}^1$.

In Karr's theorem, the random process, $\{Z(\mathbf{x}): \mathbf{x} \in \mathbb{R}^d\}$, is assumed to be a stationary Poisson process with intensity $\nu$ which is independent of the values of $Z(\mathbf{x})$. We specialize his results to our space–time random field by considering $\{Z(\mathbf{x}): \mathbf{x} \in \mathbb{R}^3\}$. Let $\mathbf{x}$ denote a location in $\mathbb{R}^3$ and $d\mathbf{x}$ denote an infinitesimally small sphere centered at $\mathbf{x}$. Define $w(\cdot)$ and $N^{(2)}(d\mathbf{x}_1, d\mathbf{x}_2)$ in the same way as in Section 3.1 with the understanding that in this section, the support of $w(\cdot)$ and $\mathbf{x}$ are defined in $\mathbb{R}^3$ rather than $\mathbb{R}^2$. The kernel estimator over $D_n$ is defined as

$$\widehat{C}_n(\mathbf{k}) = \frac{1}{\nu^2 |D_n|} \int_{D_n} \int_{D_n} w_n(\mathbf{k} - \mathbf{x}_1 + \mathbf{x}_2) Z(\mathbf{x}_1) Z(\mathbf{x}_2) N^{(2)}(d\mathbf{x}_1, d\mathbf{x}_2).$$

We define the mixing coefficients as in (1.1′) and assume the mixing condition, the boundary condition and the moment conditions (C1), (C2) and (C3′), setting $p = 0, q = 3, q_1 = 3$ and $\mathcal{I}_n = D_n$. We adopt the definition of $Q$ from (3.2) in this section and assume that $\lambda_n$ satisfies the condition

$$\lambda_n \to 0, \qquad \lambda_n^3 |D_n| \to \infty. \tag{3.3}$$

**Theorem 4.** *Let $\{Z(\mathbf{x}): \mathbf{x} \in \mathbb{R}^3\}$ be a strictly stationary mean zero spatio-temporal random field which is observed in $D_n \subset \mathbb{R}^3$ satisfying condition* (C2). *Assume that the locations $\mathbf{x}$ are generated by a homogeneous Poisson process in $\mathbb{R}^3$. Assume that $\int_{\mathbb{R}^3} C(\mathbf{k}) \, d\mathbf{k} < \infty$ and that $Q$ exists and satisfies*

$$\sup_{\mathbf{x}_1, \mathbf{x}_2} \int_{\mathbb{R}^3} |Q(\mathbf{x} + \mathbf{x}_1, \mathbf{x}, \mathbf{x}_2)| \, d\mathbf{x} < \infty.$$

*Then, $\mathrm{E}\{\widehat{C}(\mathbf{k})\} \to C(\mathbf{k})$ and $\Sigma = \lim_{n \to \infty} \lambda_n^3 |D_n| \mathrm{cov}(\widehat{\mathbf{G}}_n, \widehat{\mathbf{G}}_n)$ exists, the $(i, j)$th element of which is*

$$\frac{1}{\nu^2} \int_{\mathbb{R}^3} w^2(\mathbf{x}) \, d\mathbf{x} \, [\mathrm{E}\{Z^2(\mathbf{0}) Z^2(\mathbf{k}_i)\} I(\mathbf{k}_i = \pm \mathbf{k}_j)],$$

*where $I(\mathbf{k}_i = \pm \mathbf{k}_j) = 1$ if $\mathbf{k}_i = \pm \mathbf{k}_j$ and 0 otherwise. If we further assume that $\mathrm{E}\{Z^2(\mathbf{0}) Z^2(\mathbf{k})\} > 0$ for all $\mathbf{k} \in \Lambda$ and that conditions* (C1), (C3′) *and* (3.3) *hold, then*

$$\sqrt{\lambda_n^3 |D_n|} \{\widehat{\mathbf{G}}_n - \mathrm{E}(\widehat{\mathbf{G}}_n)\} \xrightarrow{d} \mathrm{N}_m(\mathbf{0}, \Sigma)$$

*as $n \to \infty$.*

Unlike Theorem 3, Theorem 4 combines $\mathbf{s}$ and $t$ as a single location $\mathbf{x}$. As a result, we get a more concise form for $\Sigma$ than in Theorem 3. However, we can expect there to be two separate terms for $I(\mathbf{h}_i) = I(\mathbf{h}_j)$ and $I(\mathbf{h}_i) = I(-\mathbf{h}_j)$ in the formulation of $\Sigma$ as in Theorem 3 if we consider $\mathbf{s}$ and $t$ separately.

## 4. Simulation

In order to illustrate the approach to joint normality of $\widehat{\mathbf{G}}_n$ with the derived asymptotic covariance matrix, we perform the following simulation experiment. For simplicity, we consider only the situation given in Section 2.2, noting that it corresponds to many practical applications such as the Irish wind data mentioned at the end of Section 2.2. We simulate a random field via a vector autoregressive (VAR) model with spatial structure (see, e.g., de Luna and Genton (2005)). Specifically, we assume that the fixed space $S$ is a $3 \times 3$ grid with grid interval 1 and consider the VAR(1) model

$$\mathbf{Z}_t = R\mathbf{Z}_{t-1} + \varepsilon_t, \tag{4.1}$$

where $\mathbf{Z}_t = (Z(\mathbf{s}_1, t), \ldots, Z(\mathbf{s}_9, t))^{\mathrm{T}}$ and $\varepsilon_t$ is a Gaussian multivariate white noise process with a spatial stationary and isotropic exponential correlation function given by $(\Sigma_\varepsilon)_{ij} = \exp(-\frac{\|\mathbf{s}_i - \mathbf{s}_j\|}{\phi})$, $i, j = 1, \ldots, 9$, where $\phi$ denotes the range parameter (we initially set $\phi = 1$). $R$ is a $9 \times 9$ matrix of coefficients which determines the dependency between $\mathbf{Z}_t$ and $\mathbf{Z}_{t-\tau}$ as $\mathrm{cov}(\mathbf{Z}_t, \mathbf{Z}_{t-\tau}) = R^\tau \Gamma_z(0)$, where $\mathrm{vec}\{\Gamma_z(0)\} = (I_{81} - R \otimes R)^{-1} \mathrm{vec}(\Sigma_\varepsilon)$ and $\mathrm{vec}(\cdot)$ denotes the operator vectorizing a matrix. Note that the coefficients in $R$ are determined only by $\mathbf{s}$, not by $t$. We set these coefficients as follows: for each $(\mathbf{s}_i, t)$, $i = 1, \ldots, 9$, the coefficient is 0.2 for $(\mathbf{s}_i, t-1)$, whereas it is 0.1 for $\{(\mathbf{s}_j, t-1): \|\mathbf{s}_j - \mathbf{s}_i\| = 1\}$, $1 \leq j \leq 9$, and 0 for the remaining $(\mathbf{s}, t-1)$'s.

We choose two space–time lags, $\mathbf{k}_1 = (\|\mathbf{h}\| = 1, u = 0)$ and $\mathbf{k}_2 = (\|\mathbf{h}\| = 1, u = 1)$, and set $\mathbf{\Lambda} = \{\mathbf{k}_1, \mathbf{k}_2\}$. Here, $\|\mathbf{h}\|$ denotes the Euclidean norm of $\mathbf{h}$. To explore the empirical joint distribution of $\widehat{\mathbf{G}}_n$ as $T_n = \{1, \ldots, n\}$ increases, we first set $n = 3$, simulate an $S \times T_n$ random field using (4.1) and compute $\widehat{C}_n(\mathbf{k}_1)$ and $\widehat{C}_n(\mathbf{k}_2)$ over this random field. We repeat this 5000 times to obtain a sample of $\widehat{C}_n(\mathbf{k}_1)$ and $\widehat{C}_n(\mathbf{k}_2)$ values and then compute their empirical covariance matrix. We then set $n = 10, 20, 50, 70, 100, 150, 200, 500, 1000, 5000$, respectively, and repeat the same procedure as with $n = 3$. We assess the joint normality of $\widehat{C}_n(\mathbf{k}_1)$ and $\widehat{C}_n(\mathbf{k}_2)$ for each $n$ using the multivariate measures of skewness and kurtosis introduced by Mardia (1970).

Finally, we evaluate the covariance matrix of $\widehat{\mathbf{G}}_n$. Denote by $\widehat{C}_{i,j}(u)$ the sample estimator of $C_{i,j}(u) = \mathrm{cov}\{Z(\mathbf{s}_j, t), Z(\mathbf{s}_i, t+u)\}$. Priestley ((1981), page 693) presents a formula to compute $\mathrm{cov}\{\widehat{C}_{i_1, j_1}(u_1), \widehat{C}_{i_2, j_2}(u_2)\}$ for large $n$ and Gaussian processes. The result of this is essentially an immediate form of Theorem 2 once we assume that $Z$ is a Gaussian random field and we are solely interested in $\widehat{C}_{i,j}(u)$ rather than $\widehat{C}_n(\mathbf{k})$. However, to obtain the theoretical values for the general $\widehat{C}_n(\mathbf{k})$ in Theorem 2 using Priestley's method, we need to find a link between the form in Theorem 2 and $\widehat{C}_{i,j}(u)$.

**Table 1.** Simulation results

| $n$ | $b_{1,2}$ | $b_{2,2}$ | $|T_n|\cov\{\widehat{C}_n(\mathbf{k}_1),\widehat{C}_n(\mathbf{k}_2)\}$ | $|T_n|\var\{\widehat{C}_n(\mathbf{k}_1)\}$ | $|T_n|\var\{\widehat{C}_n(\mathbf{k}_2)\}$ |
|---|---|---|---|---|---|
| 3 | 5.797 | 17.565 | 0.382 | 0.491 | 0.559 |
| 10 | 2.730 | 12.929 | 0.510 | 0.651 | 0.584 |
| 20 | 1.291 | 10.279 | 0.483 | 0.631 | 0.542 |
| 50 | 0.701 | 9.342 | 0.497 | 0.643 | 0.553 |
| 70 | 0.435 | 8.917 | 0.498 | 0.637 | 0.552 |
| 100 | 0.270 | 8.583 | 0.492 | 0.647 | 0.533 |
| 150 | 0.185 | 8.447 | 0.505 | 0.663 | 0.546 |
| 200 | 0.179 | 8.352 | 0.494 | 0.652 | 0.533 |
| 500 | 0.056 | 8.026 | 0.496 | 0.655 | 0.538 |
| 1000 | 0.043 | 8.084 | 0.511 | 0.665 | 0.554 |
| 5000 | 0.010 | 8.053 | 0.501 | 0.653 | 0.539 |
| $\infty$ | 0 | 8 | $0.497^*$ | $0.653^*$ | $0.539^*$ |

$b_{1,2}$ and $b_{2,2}$ denote multivariate measures of skewness and kurtosis, respectively.
$^*$are obtained using Lemma 1.

**Lemma 1.** *If $\mathbf{Z}_t = (Z(\mathbf{s}_1,t),\ldots,Z(\mathbf{s}_n,t))^{\mathrm{T}}$ is a Gaussian process with mean 0, $C_{i,j}(u) = \mathrm{E}\{Z(\mathbf{s}_j,t)Z(\mathbf{s}_i,t+u)\}$ and $\widehat{C}_{i,j}(u) = \frac{1}{|T_n|}\sum_t Z(\mathbf{s}_j,t)Z(\mathbf{s}_i,t+u)$, then*

$$\cov\{\widehat{C}_n(\mathbf{h}_i,u_i),\widehat{C}_n(\mathbf{h}_j,u_j)\} = \frac{1}{|S(\mathbf{h}_i)|}\frac{1}{|S(\mathbf{h}_j)|}\sum_{\mathbf{s}_k\in S(\mathbf{h}_i)}\sum_{\mathbf{s}_l\in S(\mathbf{h}_j)}\cov\{\widehat{C}_{k,k'}(u_i),\widehat{C}_{l,l'}(u_j)\},$$

*where $\mathbf{s}_{k'} = \mathbf{s}_k + \mathbf{h}_i$, $\mathbf{s}_{l'} = \mathbf{s}_l + \mathbf{h}_j$, and for large $n$,*

$$\cov\{\widehat{C}_{i_1,j_1}(s),\widehat{C}_{i_2,j_2}(u)\} \sim \frac{1}{|T_n|}\sum_{r\in\mathbb{Z}}\{C_{j_2,j_1}(r)C_{i_2,i_1}(r+u-s)+C_{i_2,j_1}(r+u)C_{j_2,i_1}(r-s)\}.$$

We calculate the theoretical values in Theorem 2 using Lemma 1 and compare them with the simulation output. The simulation output, together with the theoretical values, are summarized in Table 1. We see from this table that the multivariate measures of skewness and kurtosis approach 0 and 8 respectively as $n$ increases. This agrees with the skewness and kurtosis of a bivariate normal distribution. Additionally, all $|T_n|\cov\{\widehat{C}_n(\mathbf{k}_1),\widehat{C}_n(\mathbf{k}_2)\}$, $|T_n|\var\{\widehat{C}_n(\mathbf{k}_1)\}$ and $|T_n|\var\{\widehat{C}_n(\mathbf{k}_2)\}$ approach their corresponding theoretical values as $n$ increases. This verifies our asymptotic covariance matrix of $\widehat{\mathbf{G}}_n$. As seen in Table 1, $\cov(\widehat{\mathbf{G}}_n,\widehat{\mathbf{G}}_n)$ stabilizes at about $n=100$.

We have observed that the simulation results do not change appreciably when the number of simulations exceeds 5000, so we show only the results based on 5000 simulations. Setting the range parameter $\phi$ in the spatial correlation function to be 1.5 increases the variances and covariances, but the trends of multivariate measures of skewness and kurtosis and the convergence of $\cov(\widehat{\mathbf{G}}_n,\widehat{\mathbf{G}}_n)$ remain the same as when $\phi=1$. When $S$ is a $4\times 4$ grid, the variances and covariances become smaller than with a $3\times 3$ grid as larger grid sizes provide more data to estimate $\widehat{\mathbf{G}}_n$. However, all of our simulation experiments show a more rapid approach of $\widehat{\mathbf{G}}_n$ to joint normality as $n$ increases.

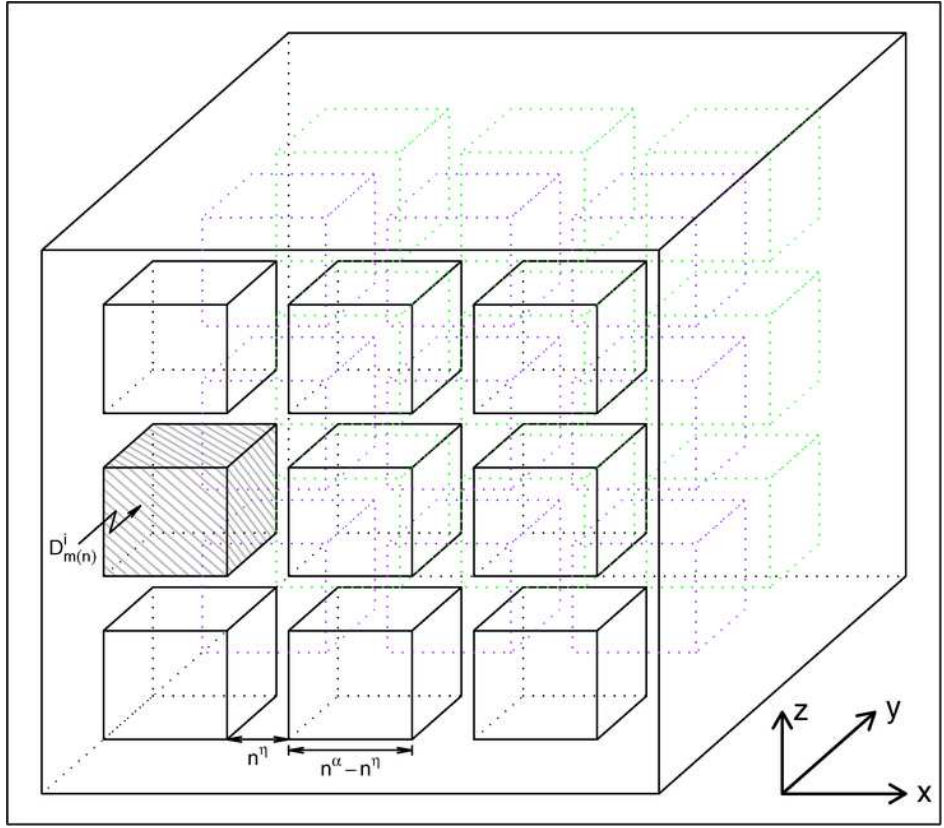

**Figure 1.** Partition of the random field for Lemma A.2.

## Appendix A: Lemmas

**Lemma A.1.** *Consider two closed and connected sets $U, V$ in $\mathbb{R}^d$ such that $|U| = |V| \leq b$ and $d(U,V) \geq r$. Let $X$ and $Y$ be measurable random variables with respect to $\mathfrak{F}(U)$ and $\mathfrak{F}(V)$ such that $|X| < C_1$ and $|Y| < C_2$. Then, $\mathrm{cov}(X,Y) \leq 4C_1 C_2 \alpha_b(r)$.*

**Proof.** The proof is completely analogous to that of Theorem 17.2.1 in Ibragimov and Linnik (1971), page 306. Note that if the variables $X$, $Y$ are complex, then separating the real and imaginary parts, we again arrive at the same expression, with 4 replaced by 16. □

**Lemma A.2.** *Let $\{Z(\mathbf{x}), \mathbf{x} \in \mathbb{R}^d\}$ be a strictly stationary mean zero random field which is observed at lattice points in $D_n \subset \mathbb{Z}^d$. Let $\pi$ be a parameter of the random field and $\widehat{\pi}$ be a consistent estimator of the form $\widehat{\pi}_n = \frac{1}{|D_n(\pi)|} \sum_{D_n(\pi)} f\{Z(\mathbf{x})\}$ satisfying $\mathrm{E}[f\{Z(\mathbf{x})\}] = \pi$, where $D_n(\pi)$ is the set of $\mathbf{x}$ being used to estimate $\pi$ and $\frac{|D_n(\pi)|}{|D_n|} \to 1$. Assume that*

$$\sum_{\mathbf{x} \in \mathbb{Z}^d} |\mathrm{cov}[f\{Z(\mathbf{0})\}, f\{Z(\mathbf{x})\}]| < \infty.$$

*Then, $\sigma^2 = \lim_{n \to \infty} |D_n| \mathrm{var}(\widehat{\pi}_n)$ exists. If we further assume conditions* (C1), (C2) *and* (C3), *setting $p = 0$, $q = d$, $\mathcal{I}_n = D_n$, $\widehat{C}_n(\mathbf{k}) = \widehat{\pi}_n$ and $C(\mathbf{k}) = \pi$, then $\sqrt{|D_n|}(\widehat{\pi}_n - \pi) \xrightarrow{d} \mathrm{N}(0, \sigma^2)$.*

**Proof.** Let $A_n = \sqrt{|D_n|}(\widehat{\pi}_n - \pi)$. Then, $\sigma^2 = \lim_{n \to \infty} \mathrm{var}(A_n) = \sum_{\mathbf{x} \in \mathbb{Z}^d} \mathrm{cov}[f\{Z(\mathbf{0})\}, f\{Z(\mathbf{x})\}]$ exists by assumption. We apply a blocking technique (e.g., Ibragimov and Linnik (1971)) in conjunction with the mixing condition to prove the normality of $\widehat{\pi}_n$.

Let $l(n) = n^\alpha$ and let $m(n) = n^\alpha - n^\eta$ for some $2d/(d+\varepsilon) < \eta < \alpha < 1$. Divide the original field $D_n$ into non-overlapping subhypercubes, $D^i_{l(n)} = l^d(n)$, $i = 1, \ldots, k_n$; within each subhypercube, further obtain $D^i_{m(n)}$ which shares the same center as $D^i_{l(n)}$. Thus, $d(D^i_{m(n)}, D^{i'}_{m(n)}) \geq n^\eta$ for $i \neq i'$. Figure 1 shows an example when $d = 3$. Let $\widehat{\pi}^i_{m(n)}$ denote the estimator obtained from $D^i_{m(n)}$. Let $a_n = \sum_{i=1}^{k_n} a^i_n / \sqrt{k_n}$, $a'_n = \sum_{i=1}^{k_n} (a^i_n)' / \sqrt{k_n}$, where $a^i_n = \sqrt{m^d(n)} \{\widehat{\pi}^i_{m(n)} - \pi\}$ and $(a^i_n)'$ have the same marginal distributions as $a^i_n$, but are independent. Let $\phi'_n(x)$ and $\phi_n(x)$ be the characteristic functions of $a'_n$ and $a_n$, respectively. The proof consists of the following three steps:

(S1) $A_n - a_n \xrightarrow{p} 0$;
(S2) $\phi'_n(x) - \phi_n(x) \to 0$;
(S3) $a'_n \xrightarrow{d} \mathrm{N}(0, \sigma^2)$.

**Proof of (S1).** Since $\mathrm{E}(A_n - a_n) = 0$, it suffices to show that $\mathrm{var}(A_n - a_n) \to 0$. Let $D^{m(n)}$ denote the union of all $D^i_{m(n)}$. Observe that

$$\begin{aligned}
a_n &= \frac{1}{\sqrt{k_n}} \sum_{i=1}^{k_n} a^i_n \\
&= \frac{\sqrt{k_n} \sqrt{m^d(n)}}{|D^{m(n)}|} \sum_{\mathbf{x} \in D^{m(n)}} [f\{Z(\mathbf{x})\} - \pi] \\
&= \sqrt{|D^{m(n)}|}(\widehat{\pi}_{D^{m(n)}} - \pi).
\end{aligned}$$

Thus, $\mathrm{var}(a_n) \to \sigma^2$ since $\frac{|D_n|}{|D^{m(n)}|} \to 1$ (demonstrated in Lemma A.3) and

$$\mathrm{cov}(A_n, a_n) = \frac{\sqrt{|D_n||D^{m(n)}|}}{|D_n||D^{m(n)}|} \sum_{\mathbf{x}_1 \in D_n} \sum_{\mathbf{x}_2 \in D^{m(n)}} \mathrm{cov}[f\{Z(\mathbf{x}_1)\}, f\{Z(\mathbf{x}_2)\}].$$

Again by $\frac{|D^{m(n)}|}{|D_n|} \to 1$, we get $\mathrm{cov}(A_n, a_n) \to \sigma^2$ and $\mathrm{var}(A_n - a_n) \to 0$.

**Proof of (S2).** We employ telescope arguments here. Let $\iota$ denote the imaginary unit. We define $U_i = \exp(\iota x \frac{a_n^i}{\sqrt{k_n}})$, $X_j = \prod_{i=1}^{j} U_i$ and $Y_j = U_{j+1}$. Analogously to Sherman and Carlstein (1994), Lemma 2, we have

$$\mathrm{cov}(X_j, Y_j) \leq 16jm^d(n)n^{-\varepsilon\eta} = 16j(n^\alpha - n^\eta)^d n^{-\varepsilon\eta}$$
$$\leq 16jn^{d\alpha - \varepsilon\eta},$$

by (C1). Thus, $|\phi_n'(x) - \phi_n(x)| \leq \sum_{j=1}^{k_n - 1} 16jn^{d\alpha - \varepsilon\eta} = \mathrm{O}(n^{2d - d\alpha - \varepsilon\eta})$. The last equality in the previous expression follows from $\mathrm{O}(k_n) = \mathrm{O}(\frac{n^d}{n^{d\alpha}}) = \mathrm{O}(n^{d(1-\alpha)})$, and $\sum_{j=1}^{k_n - 1} 16j = \mathrm{O}(k_n^2)$. Since $2d/(d+\varepsilon) < \eta < \alpha < 1$, we have $2d - d\alpha - \varepsilon\eta < 2d - d\eta - \varepsilon\eta < 0$. Then, $|\phi_n'(x) - \phi_n(x)| \to 0$.

**Proof of (S3).** Observe that $\mathrm{E}(|(a_n^i)'|^{2+\delta}) < C_\delta$ for some constant $C_\delta$. Since $(a_n^i)'$ are i.i.d., $\mathrm{var}\{\sum_{i=1}^{k_n}(a_n^i)'\} = k_n \mathrm{var}\{(a_n^i)'\}$. Defining $\sigma_n^2 = \mathrm{var}\{(a_n^i)'\}$, we have $\sigma_n^2 \to \sigma^2$ from the proof of (S1). Thus,

$$\lim_{n\to\infty} \sum_{i=1}^{k_n} \frac{\mathrm{E}(|(a_n^i)'|^{2+\delta})}{\sqrt{[\mathrm{var}\{\sum_{i=1}^{k_n}(a_n^i)'\}]^{2+\delta}}} \leq \lim_{n\to\infty} C_\delta \frac{k_n}{(k_n\sigma_n^2)^{(2+\delta)/2}} = 0.$$

Thus, applying Lyapounov's theorem, we have $\frac{1}{\sqrt{k_n}} \sum_{i=1}^{k_n} (a_n^i)' \xrightarrow{d} \mathrm{N}(0, \sigma^2)$. □

**Lemma A.3.** *Assuming that $D_n$ satisfies* (C2) *setting $\mathcal{I}_n = D_n$ and $q = d$, we have $|D^{m(n)}|/|D_n| \to 1$ as $n \to \infty$, where $D^{m(n)}$ is defined as in the proof of Lemma A.2.*

**Proof.** Introduce the following notations: $D^{l(n)}$ denotes the union of $D_{l(n)}^i$ defined in the proof of Lemma A.2; $D_{n/l(n)}$ denotes the field in $D_n$ but not in $D^{l(n)}$; $D_{l(n)/m(n)}$ denotes the field in $D^{l(n)}$ but not in $D^{m(n)}$; $k'(n)$ denotes the minimal number of extra $l^d(n)$ subhypercubes needed to cover the whole of $D_n$; $D_{n'}$ denotes the union of all extra $l^d(n)$ subhypercubes needed to cover the whole of $D_n$. Since

$$|D^{l(n)}| = k_n|D_{l(n)}^i| = k_n l^d(n) = k_n n^{d\alpha},$$
$$|D^{m(n)}| = k_n|D_{m(n)}^i| = k_n m^d(n) = k_n(n^\alpha - n^\eta)^d = k_n n^{d\alpha} + \mathrm{o}(k_n n^{d\alpha}),$$

noting that $\eta < \alpha$, we have $|D_{l(n)/m(n)}| = |D^{l(n)}| - |D^{m(n)}|$ is $o(|D^{l(n)}|)$.

We now show that $|D_{n/l(n)}|$ is $o(n^d)$. Note that $D_{n/l(n)} \subset D_{n'}$, that is, $|D_{n/l(n)}| < |D_{n'}|$. Thus, $|D_{n'}|$ being $o(n^d)$ is sufficient for $|D_{n/l(n)}|$ to be $o(n^d)$.

To show this, we split the boundary of $D_n$ into hypercubes of volume $l^{d-1}(n)$. Here, we use $L_{n,i}$ to denote the $i$th hypercube and $k_n''$ to denote the number of all the hypercubes available. $k_n''$ is $O(\frac{n^{d-1}}{l^{d-1}(n)})$ due to (C2). For each $L_{n,i}$, we may construct an $l^d(n)$ hypercube (denoted by $HC_{l(n),i}$) which fully contains $L_{n,i}$; in addition to that, we construct a $\{3l(n)\}^d$ hypercube (denoted by $HC_{3l(n),i}$) which has the same center as $HC_{l(n),i}$.

Any $l^d(n)$ hypercube that intersects $L_{n,i}$ (and thus intersects $HC_{l(n),i}$) will be fully contained in $HC_{3l(n),i}$. Since every subhypercube in $D_{n'}$ intersects the boundary of $D_n$ and thus at least one of these $L_{n,i}$'s, it will be fully contained in one of the $HC_{3l(n),i}$'s. These $l^d(n)$ hypercubes in $D_{n'}$ do not intersect each other, except at the boundary. The maximum number of such hypercubes that intersect $L_{n,i}$ cannot be larger than $3^d$. Since the size of $k_n''$ is $O(\frac{n^{d-1}}{l^{d-1}(n)})$, we conclude that $k_n'$ is no larger than $O(\frac{n^{d-1}}{l^{d-1}(n)})$ due to $k_n' < 3^d k_n''$. Thus, $|D_{n/l(n)}|$ is no larger than $O(\frac{n^{d-1}}{l^{d-1}(n)} l^d(n)) = O(n^{d-1} l(n)) = o(n^d)$. Then,

$$|D_n| - |D^{m(n)}| = \{|D_n| - |D^{l(n)}|\} + \{|D^{l(n)}| - |D^{m(n)}|\}$$
$$= |D_{n/l(n)}| + |D_{l(n)/m(n)}| = o(|D_n|).$$

Thus, $|D^{m(n)}|/|D_n| \to 1$ as $n \to \infty$. □

**Lemma A.4.** *Let $\{Z(\mathbf{x}), \mathbf{x} \in \mathbb{R}^d\}$ be a strictly stationary mean zero random field which is observed in $D_n = S_n \times T_n$, where $S_n \subset \mathbb{R}^p$, $T_n \subset \mathbb{Z}^q$ and $p + q = d$. Assume that the locations of observations in $T_n$ are regularly spaced and that the locations in $S_n$ are generated by a homogeneous Poisson process. Let $|S_n|$ denote the volume of $S_n$ and $|T_n|$ the cardinality of $T_n$. Let*

$$\widehat{C}_n(\mathbf{k}) = \frac{1}{\nu^2 |T_n| |S_n|} \sum_{\mathbf{x}_q \in T_n} \int_{S_n} \int_{S_n} w_n(\mathbf{k}_p - \mathbf{x}_{1p} + \mathbf{x}_{2p}) Z(\mathbf{x}_{1p}, \mathbf{x}_q) Z(\mathbf{x}_{2p}, \mathbf{x}_q + \mathbf{k}_q)$$
$$\times N^{(2)}(d\mathbf{x}_{1p}, d\mathbf{x}_{2p}),$$

*where $\mathbf{k} = (\mathbf{k}_p^T, \mathbf{k}_q^T)^T$ and $\mathbf{x} = (\mathbf{x}_p^T, \mathbf{x}_q^T)^T$. Assume that*

$$\sup_{\mathbf{k}_q} \int_{\mathbb{R}^p} C(\mathbf{k}_p, \mathbf{k}_q) \, d\mathbf{k}_p < \infty,$$

$$\sup_{\mathbf{x}_p, \mathbf{x}_q, \mathbf{x}_{2q}} \sum_{\mathbf{x}_{1q} \in \mathbb{Z}^q} \mathrm{E}\{Z(\mathbf{0}, \mathbf{0}) Z(\mathbf{x}_p, \mathbf{x}_q) Z(\mathbf{0}, \mathbf{x}_{1q}) Z(\mathbf{x}_p, \mathbf{x}_{1q} + \mathbf{x}_{2q})\} < \infty,$$

$$\sup_{\mathbf{x}_{1p}, \mathbf{x}_{2p}} \sup_{\mathbf{x}_q, \mathbf{x}_{1q}, \mathbf{x}_{2q}} \int_{\mathbf{x}_p \in \mathbb{R}^p} |Q\{(\mathbf{x}_{1p}, \mathbf{x}_{1q}), (\mathbf{x}_p, \mathbf{x}_q), (\mathbf{x}_p + \mathbf{x}_{2p}), (\mathbf{x}_q + \mathbf{x}_{2q})\}| \, d\mathbf{x}_{2p} < \infty,$$

where $w_n(\mathbf{x}) = \frac{1}{\lambda_n^p} w(\frac{\mathbf{x}}{\lambda_n})$ and $\lambda_n$ is a sequence of positive constants satisfying $\lambda_n \to 0$ and $\lambda_n^p |S_n| \to \infty$. $Q$ is defined as in (3.2). Then, $\mathrm{E}\{\widehat{C}_n(\mathbf{k})\} \to C(\mathbf{k})$ and $\Sigma = \lim_{n \to \infty} |T_n| \times |S_n| \lambda_n^p \mathrm{cov}(\widehat{\mathbf{G}}_n, \widehat{\mathbf{G}}_n)$ exists, the $(i,j)$th element of which is

$$\frac{1}{\nu^2} \int w^2(\mathbf{y}) \, \mathrm{d}\mathbf{y} \sum_{\mathbf{x}_q \in \mathbb{Z}^q} \mathrm{E}[Z(\mathbf{0},\mathbf{0})Z(\mathbf{k}_{ip}, \mathbf{k}_{iq}) \{Z(\mathbf{0},\mathbf{x}_q)Z(\mathbf{k}_{ip}, \mathbf{x}_q + \mathbf{k}_{jq}) I(\mathbf{k}_{ip} = \mathbf{k}_{jp})$$
$$+ Z(\mathbf{k}_{ip}, \mathbf{x}_q) Z(\mathbf{0}, \mathbf{x}_q + \mathbf{k}_{jq}) I(\mathbf{k}_{ip} = -\mathbf{k}_{jp})\}],$$

where $I(\mathbf{k}_{ip} = \pm \mathbf{k}_{jp}) = 1$ if $\mathbf{k}_{ip} = \pm \mathbf{k}_{jp}$ and 0 otherwise. If we further assume that $\Sigma$ is positive definite and that conditions (C1), (C2) and (C3$'$) hold, setting $p = 0$, $q = d$, $q_1 = p$, $\mathcal{I}_n = D_n$ and replacing $C(\mathbf{k})$ with $\mathrm{E}\{\widehat{C}_n(\mathbf{k})\}$, then

$$\sqrt{|T_n||S_n|\lambda_n^p} \{\widehat{\mathbf{G}}_n - \mathrm{E}(\widehat{\mathbf{G}}_n)\} \xrightarrow{d} \mathrm{N}_m(\mathbf{0}, \Sigma).$$

**Proof.**

$\mathrm{E}\{\widehat{C}_n(\mathbf{k})\}$
$$= \frac{1}{\nu^2 |T_n||S_n|} \sum_{\mathbf{x}_q \in T_n} \mathrm{E}\left\{ \iint_{S_n} w_n(\mathbf{k}_p + \mathbf{x}_{1p} - \mathbf{x}_{2p}) Z(\mathbf{x}_{1p}, \mathbf{x}_q) Z(\mathbf{x}_{2p}, \mathbf{x}_q + \mathbf{k}_q) N^{(2)}(\mathrm{d}\mathbf{x}_{1p}, \mathrm{d}\mathbf{x}_{2p}) \right\}$$
$$= \frac{1}{|T_n|} \sum_{\mathbf{x}_q \in T_n} \frac{1}{|S_n|} \iint_{S_n} \frac{1}{\lambda_n^p} w\left(\frac{\mathbf{k}_p + \mathbf{x}_{1p} - \mathbf{x}_{2p}}{\lambda_n}\right) C(\mathbf{x}_{2p} - \mathbf{x}_{1p}, \mathbf{k}_q) \, \mathrm{d}\mathbf{x}_{1p} \, \mathrm{d}\mathbf{x}_{2p}$$
$$= \frac{1}{|T_n|} \sum_{\mathbf{x}_q \in T_n} \int_{S_n - S_n} w(\mathbf{y}) C(\mathbf{k}_p - \lambda_n \mathbf{y}, \mathbf{k}_q) \frac{|S_n \cap (S_n - \mathbf{k}_p + \lambda_n \mathbf{y})|}{|S_n|} \, \mathrm{d}\mathbf{y},$$

where $S_n - S_n$ denotes all pairwise differences between the two sets. Since

$$\int_{S_n - S_n} w(\mathbf{y}) C(\mathbf{k}_p - \lambda_n \mathbf{y}, \mathbf{k}_q) \frac{|S_n \cap (S_n - \mathbf{k}_p + \lambda_n \mathbf{y})|}{|S_n|} \, \mathrm{d}\mathbf{y}$$
$$\longrightarrow C(\mathbf{k}_p, \mathbf{k}_q) \int_{\mathbb{R}^q} w(\mathbf{y}) \, \mathrm{d}\mathbf{y} = C(\mathbf{k}_p, \mathbf{k}_q),$$

we have $\mathrm{E}\{\widehat{C}_n(\mathbf{k})\} \longrightarrow C(\mathbf{k})$.

The derivation for the variance is analogous to Karr (1986). Specifically, consider two spatial locations $\mathbf{k} = (\mathbf{k}_p^{\mathrm{T}}, \mathbf{k}_q^{\mathrm{T}})^{\mathrm{T}}$ and $\mathbf{k}' = (\mathbf{k}_p'^{\mathrm{T}}, \mathbf{k}_q'^{\mathrm{T}})^{\mathrm{T}}$:

$\mathrm{E}\{\widehat{C}_n(\mathbf{k})\widehat{C}_n(\mathbf{k}')\}$
$$= \frac{1}{\nu^4 |T_n|^2} \sum_{\mathbf{x}_q \in T_n} \sum_{\mathbf{x}_q' \in T_n} \frac{1}{|S_n|^2} \iiiint_{S_n} w_n(\mathbf{k}_p + \mathbf{x}_{1p} - \mathbf{x}_{2p}) w_n(\mathbf{k}_p' + \mathbf{x}_{1p}' - \mathbf{x}_{2p}')$$

$$\times \mathrm{E}\{Z(\mathbf{x}_{1p},\mathbf{x}_q)Z(\mathbf{x}_{2p},\mathbf{x}_q+\mathbf{k}_q)Z(\mathbf{x}'_{1p},\mathbf{x}'_q) \quad \text{(A.1)}$$
$$\times Z(\mathbf{x}'_{2p},\mathbf{x}'_q+\mathbf{k}'_q)\}$$
$$\times \mathrm{E}\{N^{(2)}(\mathrm{d}\mathbf{x}_{1p},\mathrm{d}\mathbf{x}_{2p})N^{(2)}(\mathrm{d}\mathbf{x}'_{1p},\mathrm{d}\mathbf{x}'_{2p})\}.$$

Observe that

$$\begin{aligned}
&\mathrm{E}\{Z(\mathbf{x}_{1p},\mathbf{x}_q)Z(\mathbf{x}_{2p},\mathbf{x}_q+\mathbf{k}_q)Z(\mathbf{x}'_{1p},\mathbf{x}'_q)Z(\mathbf{x}'_{2p},\mathbf{x}'_q+\mathbf{k}'_q)\} \\
&= Q\{(\mathbf{x}_{2p}-\mathbf{x}_{1p},\mathbf{k}_q),(\mathbf{x}'_{1p}-\mathbf{x}_{1p},\mathbf{x}'_q-\mathbf{x}_q),(\mathbf{x}'_{2p}-\mathbf{x}_{1p},\mathbf{x}'_q-\mathbf{x}_q+\mathbf{k}'_q)\} \\
&\quad + C(\mathbf{x}_{2p}-\mathbf{x}_{1p},\mathbf{k}_q)C(\mathbf{x}'_{2p}-\mathbf{x}'_{1p},\mathbf{k}'_q) \\
&\quad + C(\mathbf{x}'_{1p}-\mathbf{x}_{1p},\mathbf{x}'_q-\mathbf{x}_q)C(\mathbf{x}'_{2p}-\mathbf{x}_{2p},\mathbf{x}'_q-\mathbf{x}_q+\mathbf{k}'_q-\mathbf{k}_q) \\
&\quad + C(\mathbf{x}'_{2p}-\mathbf{x}_{1p},\mathbf{x}'_q-\mathbf{x}_q+\mathbf{k}'_q)C(\mathbf{x}_{2p}-\mathbf{x}'_{1p},\mathbf{x}_q-\mathbf{x}'_q+\mathbf{k}_q).
\end{aligned} \quad \text{(A.2)}$$

In addition,

$$\begin{aligned}
&\mathrm{E}\{N^{(2)}(\mathrm{d}\mathbf{x}_{1p},\mathrm{d}\mathbf{x}_{2p})N^{(2)}(\mathrm{d}\mathbf{x}'_{1p},\mathrm{d}\mathbf{x}'_{2p})\} \\
&= \nu^4\,\mathrm{d}\mathbf{x}_{1p}\,\mathrm{d}\mathbf{x}_{2p}\,\mathrm{d}\mathbf{x}'_{1p}\,\mathrm{d}\mathbf{x}'_{2p} \\
&\quad + \nu^3\,\mathrm{d}\mathbf{x}_{1p}\,\mathrm{d}\mathbf{x}_{2p}\,\varepsilon_{\mathbf{x}_{1p}}(\mathrm{d}\mathbf{x}'_{1p})\,\mathrm{d}\mathbf{x}'_{2p} + \nu^3\,\mathrm{d}\mathbf{x}_{1p}\,\mathrm{d}\mathbf{x}_{2p}\,\mathrm{d}\mathbf{x}'_{1p}\,\varepsilon_{\mathbf{x}_{1p}}(\mathrm{d}\mathbf{x}'_{2p}) \\
&\quad + \nu^3\,\mathrm{d}\mathbf{x}_{1p}\,\mathrm{d}\mathbf{x}_{2p}\,\varepsilon_{\mathbf{x}_{2p}}(\mathrm{d}\mathbf{x}'_{1p})\,\mathrm{d}\mathbf{x}'_{2p} + \nu^3\,\mathrm{d}\mathbf{x}_{1p}\,\mathrm{d}\mathbf{x}_{2p}\,\mathrm{d}\mathbf{x}'_{1p}\,\varepsilon_{\mathbf{x}_{2p}}(\mathrm{d}\mathbf{x}'_{2p}) \\
&\quad + \nu^2\,\mathrm{d}\mathbf{x}_{1p}\,\mathrm{d}\mathbf{x}_{2p}\,\varepsilon_{\mathbf{x}_{1p}}(\mathrm{d}\mathbf{x}'_{2p})\varepsilon_{\mathbf{x}_{1p}}(\mathrm{d}\mathbf{x}'_{2p}) + \nu^2\,\mathrm{d}\mathbf{x}_{1p}\,\mathrm{d}\mathbf{x}_{2p}\,\varepsilon_{\mathbf{x}_{2p}}(\mathrm{d}\mathbf{x}'_{1p})\varepsilon_{\mathbf{x}_{2p}}(\mathrm{d}\mathbf{x}'_{1p}),
\end{aligned} \quad \text{(A.3)}$$

where $\varepsilon_{\mathbf{x}}(\cdot)$ denotes point measure.

Substituting these two expansions into the main formula produces a lengthy expression that we do not reproduce in its entirety. We only show some representative terms, rather than the total 28 terms. For simplicity, we first look at only the integral part conditional on $\mathbf{x}_q$ and $\mathbf{x}'_q$. Expanding (A.1), first by (A.2) and then by (A.3), the first term is

$$\frac{1}{\nu^4|S_n|^2}\iiiint\limits_{S_n} w_n(\mathbf{k}_p+\mathbf{x}_{1p}-\mathbf{x}_{2p})w_n(\mathbf{k}'_p+\mathbf{x}'_{1p}-\mathbf{x}'_{2p})$$
$$\times Q\{(\mathbf{x}_{2p}-\mathbf{x}_{1p},\mathbf{k}_q),(\mathbf{x}'_{1p}-\mathbf{x}_{1p},\mathbf{x}'_q-\mathbf{x}_q),$$
$$(\mathbf{x}'_{2p}-\mathbf{x}_{1p},\mathbf{x}'_q-\mathbf{x}_q+\mathbf{k}'_q)\}\nu^4\,\mathrm{d}\mathbf{x}_{1p}\,\mathrm{d}\mathbf{x}_{2p}\,\mathrm{d}\mathbf{x}'_{1p}\,\mathrm{d}\mathbf{x}'_{2p}$$
$$\leq \iiint\limits_{S_n-S_n} w_n(\mathbf{k}_p-\mathbf{v}_1)w_n(\mathbf{k}'_p+\mathbf{v}_2-\mathbf{v}_3)$$
$$\times |Q\{(\mathbf{v}_1,\mathbf{k}_q),(\mathbf{v}_2,\mathbf{x}'_q-\mathbf{x}_q),(\mathbf{v}_3,\mathbf{x}'_q-\mathbf{x}_q+\mathbf{k}'_q)\}|\frac{1}{|S_n|}\,\mathrm{d}\mathbf{v}_1\,\mathrm{d}\mathbf{v}_2\,\mathrm{d}\mathbf{v}_3$$
$$\leq C_1\iint\limits_{\mathbb{R}^p} w_n(\mathbf{k}_p-\mathbf{v}_1)w_n(\mathbf{k}'_p-\mathbf{v}_4)\frac{1}{|S_n|}\,\mathrm{d}\mathbf{v}_1\,\mathrm{d}\mathbf{v}_4$$

$$= O\left(\frac{1}{|S_n|}\right).$$

The second term is

$$\frac{1}{\nu^4 |S_n|^2} \iiiint_{S_n} w_n(\mathbf{k}_p + \mathbf{x}_{1p} - \mathbf{x}_{2p}) w_n(\mathbf{k}'_p + \mathbf{x}'_{1p} - \mathbf{x}'_{2p})$$
$$\times C(\mathbf{x}_{2p} - \mathbf{x}_{1p}, \mathbf{k}_q) C(\mathbf{x}'_{2p} - \mathbf{x}'_{1p}, \mathbf{k}'_q) \nu^4 \, d\mathbf{x}_{1p} \, d\mathbf{x}_{2p} \, d\mathbf{x}'_{1p} \, d\mathbf{x}'_{2p}$$
$$= \int_{S_n - S_n} w_n(\mathbf{k}_p - \mathbf{v}_1) C(\mathbf{v}_1, \mathbf{k}_q) \frac{|S_n \cap (S_n - \mathbf{v}_1)|}{|S_n|} \, d\mathbf{v}_1$$
$$\times \int_{S_n - S_n} w_n(\mathbf{k}'_p - \mathbf{v}_2) C(\mathbf{v}_2, \mathbf{k}'_q) \frac{|S_n \cap (S_n - \mathbf{v}_2)|}{|S_n|} \, d\mathbf{v}_2$$
$$\to C(\mathbf{k}) \int_{\mathbb{R}^p} w(\mathbf{y}) \, d\mathbf{y} \times C(\mathbf{k}') \int_{\mathbb{R}^p} w(\mathbf{y}) \, d\mathbf{y}$$
$$= C(\mathbf{k})|_{\mathbf{x}_q} \times C(\mathbf{k}')|_{\mathbf{x}'_q},$$

where $C(\mathbf{k})|_{\mathbf{x}_q}$ denotes $C(\mathbf{k})$ conditional on $\mathbf{x}_q$. Note that

$$\frac{1}{|T_n|^2} \sum_{\mathbf{x}_q \in T_n} \sum_{\mathbf{x}'_q \in T_n} C(\mathbf{k})|_{\mathbf{x}_q} \times C(\mathbf{k}')|_{\mathbf{x}'_q} \to C(\mathbf{k}) C(\mathbf{k}').$$

Most of the terms are smaller than or the same as $O(\frac{1}{|S_n|})$. Only 8 among the 28 terms are larger than $O(\frac{1}{|S_n|})$ and thus considered as dominant contributions. Combining those 8 terms, we obtain

$$|S_n| \lambda_n^p \operatorname{cov}\{C(\mathbf{k})|_{\mathbf{x}_q}, C(\mathbf{k}')|_{\mathbf{x}'_q}\}$$
$$\longrightarrow \frac{1}{\nu^2} \int_{\mathbb{R}^p} w^2(\mathbf{y}) \, d\mathbf{y} \, E[Z(\mathbf{0}, \mathbf{0}) Z(\mathbf{k}_p, \mathbf{k}_q)$$
$$\times \{Z(\mathbf{0}, \mathbf{x}'_q - \mathbf{x}_q) Z(\mathbf{k}_p, \mathbf{x}'_q - \mathbf{x}_q + \mathbf{k}'_q) I(\mathbf{k}_p = \mathbf{k}'_p)$$
$$+ Z(\mathbf{k}_p, \mathbf{x}'_q - \mathbf{x}_q) Z(\mathbf{0}, \mathbf{x}'_q - \mathbf{x}_q + \mathbf{k}'_q) I(\mathbf{k}_p = -\mathbf{k}'_p)\}].$$

Let

$$X_n := \frac{1}{|T_n|^2} \sum_{\mathbf{x}_q \in T_n} \sum_{\mathbf{x}'_q \in T_n} E[Z(\mathbf{0}, \mathbf{0}) Z(\mathbf{k}_p, \mathbf{k}_q)$$
$$\times \{Z(\mathbf{0}, \mathbf{x}'_q - \mathbf{x}_q) Z(\mathbf{k}_p, \mathbf{x}'_q - \mathbf{x}_q + \mathbf{k}'_q) I(\mathbf{k}_p = \mathbf{k}'_p)$$
$$+ Z(\mathbf{k}_p, \mathbf{x}'_q - \mathbf{x}_q) Z(\mathbf{0}, \mathbf{x}'_q - \mathbf{x}_q + \mathbf{k}'_q) I(\mathbf{k}_p = -\mathbf{k}'_p)\}]$$
$$= E\bigg[Z(\mathbf{0}, \mathbf{0}) Z(\mathbf{k}_p, \mathbf{k}_q) \sum_{T_n - T_n} \{Z(\mathbf{0}, \mathbf{v}) Z(\mathbf{k}_p, \mathbf{v} + \mathbf{k}'_q) I(\mathbf{k}_p = \mathbf{k}'_p)$$

$$+ Z(\mathbf{k}_p, \mathbf{v}) Z(\mathbf{0}, \mathbf{v} + \mathbf{k}'_q) I(\mathbf{k}_p = -\mathbf{k}'_p)\} \frac{T_n \cap (T_n - \mathbf{v})}{|T_n|^2}\bigg].$$

Applying Kronecker's lemma, we have

$$|T_n| X_n \longrightarrow \sum_{\mathbf{x}_q \in \mathbb{Z}^q} \mathrm{E}[Z(\mathbf{0}) Z(\mathbf{k}) \{ Z(\mathbf{0}, \mathbf{x}_q) Z(\mathbf{k}_p, \mathbf{x}_q + \mathbf{k}'_q) I(\mathbf{k}_p = \mathbf{k}'_p)$$
$$+ Z(\mathbf{k}_p, \mathbf{x}_q) Z(\mathbf{0}, \mathbf{x}_q + \mathbf{k}'_q) I(\mathbf{k}_p = -\mathbf{k}'_p)\}].$$

Thus,

$$|T_n||S_n|\lambda_n^p \operatorname{cov}\{\widehat{C}_n(\mathbf{k}), \widehat{C}_n(\mathbf{k}')\}$$
$$\longrightarrow \sum_{\mathbf{x}_q \in \mathbb{Z}^q} \mathrm{E}[Z(\mathbf{0}) Z(\mathbf{k}) \{ Z(\mathbf{0}, \mathbf{x}_q) Z(\mathbf{k}_p, \mathbf{x}_q + \mathbf{k}'_q) I(\mathbf{k}_p = \mathbf{k}'_p)$$
$$+ Z(\mathbf{k}_p, \mathbf{x}_q) Z(\mathbf{0}, \mathbf{x}_q + \mathbf{k}'_q) I(\mathbf{k}_p = -\mathbf{k}'_p)\}]$$
$$\times \frac{1}{\nu^2} \int_{\mathbb{R}^p} w^2(\mathbf{y}) \, \mathrm{d}\mathbf{y}.$$

There are two different terms for $I(\mathbf{k}_p = \mathbf{k}'_p)$ and $I(\mathbf{k}_p = -\mathbf{k}'_p)$ because $C(\mathbf{k}_p, \mathbf{k}_q)$ does not equal $C(-\mathbf{k}_p, \mathbf{k}_q)$ in general.

The proof of normality again uses three steps, as in Lemma A.2. The proof of (S1) is analogous to that in Lemma A.2, but uses a kernel estimator and an additional requirement for $\alpha$, $\lambda_n^p n^{\alpha p} \to \infty$. We follow Politis, Paparoditis and Romano (1998) to prove (S2). Define $X_i$ and $Y_i$ as in Lemma A.2 and

$$\mathrm{E}^N(X_i) = \mathrm{E}(X_i | N), \qquad \mathrm{E}^N(Y_i) = \mathrm{E}(Y_i | N), \qquad \operatorname{cov}^N(X_i, Y_i) = \operatorname{cov}(X_i, Y_i | N).$$

Then,

$$\operatorname{cov}(X_i, Y_i) = \mathrm{E}\{\operatorname{cov}^N(X_i, Y_i)\} + \operatorname{cov}\{\mathrm{E}^N(X_i), \mathrm{E}^N(Y_i)\}.$$

Since $N$ is a homogeneous Poisson process and $X_i$, $Y_i$ are random variables defined on two disjoint random fields, we have $\operatorname{cov}\{\mathrm{E}^N(X_i), \mathrm{E}^N(Y_i)\} = 0$. For given $N$, $X_i | N$ is measurable with respect to $\mathfrak{F}(\bigcup_{j=1}^{i} D_{m(n)}^j)$ and $Y_i | N$ is measurable with respect to $\mathfrak{F}(D_{m(n)}^{i+1})$. We then have $|\phi'_n(x) - \phi_n(x)| \leq \operatorname{cov}(X_i, Y_i | N)$. The rest of the proofs of (S2) and (S3) are completely analogous to Lemma A.2. The proof of the joint normality follows from the Cramér–Wold device. □

**Lemma A.5.** *Let $\{Z(\mathbf{x}), \mathbf{x} \in \mathbb{R}^d\}$ be a strictly stationary random field with mean $\mu$ which is observed in $D_n \subset \mathbb{R}^d$ satisfying conditions* (C1) *and* (C2), *setting $q = d$ and $\mathcal{I}_n = D_n$. Let $\overline{Z}_n = \frac{1}{|D_n|} \sum_{D_n} Z(\mathbf{x})$. If*

$$\sum_{\mathbf{x} \in \mathbb{R}^d} |\operatorname{cov}\{Z(\mathbf{0}), Z(\mathbf{x})\}| < \infty \qquad \textit{for all finite } \mathbf{x}$$

*and*

$$\sup_n \mathrm{E}|\sqrt{|D_n|}(\overline{Z}-\mu)|^{2+\delta} \leq C_\delta, \qquad \text{for some } \delta > 0, C_\delta < \infty,$$

*then* $\sigma_{\overline{Z}}^2 = \lim_{n\to\infty} |D_n| \mathrm{var}(\overline{Z}_n)$ *exists and* $\sqrt{|D_n|}(\overline{Z}_n - \mu) \xrightarrow{d} \mathrm{N}(0, \sigma_{\overline{Z}}^2)$ *as* $n \to \infty$.

**Proof.** The proof follows directly from Lemma A.2. □

It is straightforward to prove that the results in Lemma A.5 hold for all of the situations we have considered in this article, given the appropriate mixing condition. We have proven that $\widehat{\mathbf{G}}_n = \{\widehat{C}_n(\mathbf{h}, u) : (\mathbf{h}, u) \in \mathbf{\Lambda}\}$ is asymptotically jointly normal, that is, $\sqrt{|D_n|} \times (\widehat{\mathbf{G}}_n - \widehat{\mathbf{G}}) \xrightarrow{d} \mathrm{N}_m(\mathbf{0}, \Sigma)$ under appropriate assumptions. In this formulation of $\widehat{\mathbf{G}}_n$, we have assumed the mean of $Z(\mathbf{s}, t)$, $\mu(Z)$, is known and, without loss of generality, assumed it to equal 0. However, in most cases, the mean of the random process is unknown and needs to be estimated.

**Lemma A.6.** *Let* $\widehat{\mu}_n := \overline{Z}_n = \frac{1}{|D_n|}\sum_{\mathbf{x}\in D_n} Z(\mathbf{x})$ *and let* $\widehat{\mathbf{G}}_n^* = \{\widehat{C}_n^*(\mathbf{k}) : \mathbf{k} \in \mathbf{\Lambda}\}$, *where* $\widehat{C}_n^*(\mathbf{k}) = \frac{1}{|D_n(\mathbf{k})|}\sum_{\mathbf{x}\in D_n(\mathbf{k})} \{Z(\mathbf{x}) - \widehat{\mu}_n\}\{Z(\mathbf{x}+\mathbf{k}) - \widehat{\mu}_n\}$. *The vector of sample covariances* $\widehat{\mathbf{G}}_n^*$ *then has the same asymptotic properties as* $\widehat{\mathbf{G}}_n$, *that is,* $\sqrt{|D_n|}(\widehat{\mathbf{G}}_n^* - \mathbf{G}) \xrightarrow{d} \mathrm{N}_m(\mathbf{0}, \Sigma)$.

**Proof.** Assume that $\mu \neq 0$. We now prove that $\sqrt{|D_n|}\{\widehat{C}_n(\mathbf{k}) - \widehat{C}_n^*(\mathbf{k})\} = \mathrm{o}_p(1)$ as $n \to \infty$. Following Brockwell and Davis (1991), page 229, we have

$$\sqrt{|D_n|}\{\widehat{C}_n(\mathbf{k}) - \widehat{C}_n^*(\mathbf{k})\} = \sqrt{|D_n|}\frac{1}{|D_n(\mathbf{k})|}\sum_{\mathbf{x}\in D_n(\mathbf{k})}(\overline{Z}_n - \mu)\{Z(\mathbf{x}) + Z(\mathbf{x}+\mathbf{k})\} + (\mu^2 - \overline{Z}_n^2)$$

$$= \sqrt{|D_n|}(\overline{Z}_n - \mu)\frac{1}{|D_n(\mathbf{k})|}\sum_{\mathbf{x}\in D_n(\mathbf{k})}\{Z(\mathbf{x}) + Z(\mathbf{x}+\mathbf{k}) - \mu - \overline{Z}_n\}$$

$$= \sqrt{|D_n|}(\overline{Z}_n - \mu)(\overline{Z}_n - \mu).$$

By Lemma A.5, we know that $\sqrt{|D_n|}(\overline{Z}_n - \mu) \xrightarrow{d} \mathrm{N}(0, \sigma_{\overline{Z}}^2)$, which implies that $\sqrt{|D_n|}(\overline{Z}_n - \mu)$ is $\mathrm{O}_p(1)$ and $\overline{Z}_n - \mu \xrightarrow{p} 0$.

From these observations, we conclude that $\sqrt{|D_n|}(\widehat{\mathbf{G}}_n - \widehat{\mathbf{G}}_n^*) = \mathrm{o}_p(1)$. Thus, $\sqrt{|D_n|}(\widehat{\mathbf{G}}_n^* - \mathbf{G}) \xrightarrow{d} \mathrm{N}_m(\mathbf{0}, \Sigma)$. □

**Proof of Lemma 1.** When $Z_t$ is Gaussian, the fourth order cumulant function $Q = 0$ (Isserlis (1918)). Hence,

$$\mathrm{cov}\{\widehat{C}_{i_1,j_1}(s), \widehat{C}_{i_2,j_2}(u)\}$$
$$= \mathrm{E}\{\widehat{C}_{i_1,j_1}(s)\widehat{C}_{i_2,j_2}(u)\} - \mathrm{E}\{\widehat{C}_{i_1,j_1}(s)\}\mathrm{E}\{\widehat{C}_{i_2,j_2}(u)\}$$

$$= \frac{1}{|T_n|^2} \sum_{t_1} \sum_{t_2} \mathrm{E}(Z_{j_1,t_1} Z_{i_1,t_1+s} Z_{j_2,t_2} Z_{i_2,t_2+u}) - C_{i_1,j_1}(s) C_{i_2,j_2}(u)$$

$$= \frac{1}{|T_n|^2} \sum_{t_1} \sum_{t_2} \{C_{i_1,j_1}(s) C_{i_2,j_2}(u) + C_{j_2,j_1}(t_2 - t_1) C_{i_2,i_1}(t_2 - t_1 + u - s)$$

$$+ C_{i_2,j_1}(t_2 - t_1 + u) C_{j_2,i_1}(t_2 - t_1 - s)\} - C_{i_1,j_1}(s) C_{i_2,j_2}(u) + Q.$$

So, $|T_n| \mathrm{cov}\{\widehat{C}_{i_1,j_1}(s), \widehat{C}_{i_2,j_2}(u)\} \to \sum_{r \in \mathbb{Z}} \{C_{j_2,j_1}(r) C_{i_2,i_1}(r+u-s) + C_{i_2,j_1}(r+u) C_{j_2,i_1}(r-s)\}$ as $n \to \infty$. Let $T_n(u) = \{t : t \in T_n, t+u \in T_n\}$. We then have

$$\mathrm{cov}\{\widehat{C}_n(\mathbf{h}_i, u_i), \widehat{C}_n(\mathbf{h}_j, u_j)\}$$

$$= \frac{1}{|S(\mathbf{h}_i)|} \frac{1}{|S(\mathbf{h}_j)|} \sum_{\mathbf{s}_k \in S(\mathbf{h}_i)} \sum_{\mathbf{s}_l \in S(\mathbf{h}_j)} \mathrm{cov}\bigg\{ \frac{1}{|T_n|} \sum_{t_k \in T_n(u_i)} Z(\mathbf{s}_k, t_k) Z(\mathbf{s}_k + \mathbf{h}_i, t_k + u_i),$$

$$\times \frac{1}{|T_n|} \sum_{t_l \in T_n(u_j)} Z(\mathbf{s}_l, t_l) Z(\mathbf{s}_l + \mathbf{h}_j, t_l + u_j) \bigg\}$$

$$= \frac{1}{|S(\mathbf{h}_i)|} \frac{1}{|S(\mathbf{h}_j)|} \sum_{\mathbf{s}_k \in S(\mathbf{h}_i)} \sum_{\mathbf{s}_l \in S(\mathbf{h}_j)} \mathrm{cov}\{\widehat{C}_{k,k'}(u_i), \widehat{C}_{l,l'}(u_j)\}. \qquad \square$$

## Appendix B: Proofs of the theorems

**Proof of Theorem 1.** Consider the covariance term, $\mathrm{cov}\{\widehat{C}_n(\mathbf{h}_i, u_i), \widehat{C}_n(\mathbf{h}_j, u_j)\}$, where $(\mathbf{h}_i, u_i), (\mathbf{h}_j, u_j) \in \Lambda$, and let $D_n^i = D_n(\mathbf{h}_i, u_i)$. Then,

$$\frac{1}{|D_n^i||D_n^j|} \sum_{D_n^i} \sum_{D_n^j} \mathrm{cov}\{Z(\mathbf{s}_1, t_1) Z(\mathbf{s}_1 + \mathbf{h}_i, t_1 + u_i), Z(\mathbf{s}_2, t_2) Z(\mathbf{s}_2 + \mathbf{h}_j, t_1 + u_j)\}$$

$$= \frac{1}{|D_n^i||D_n^j|} \sum_{(\mathbf{s},t) \in D_n^j - D_n^i} \sum_{D_n^i \cap (D_n^j - (\mathbf{s},t))} \mathrm{cov}\{Z(\mathbf{0},0) Z(\mathbf{h}_i, u_i), Z(\mathbf{s},t) Z(\mathbf{s} + \mathbf{h}_j, t + u_j)\}$$

$$= \sum_{D_n^j - D_n^i} \mathrm{cov}\{Z(\mathbf{0},0) Z(\mathbf{h}_i, u_i), Z(\mathbf{s},t) Z(\mathbf{s} + \mathbf{h}_j, t + u_j)\} \frac{|D_n^i \cap \{D_n^j - (\mathbf{s},t)\}|}{|D_n^i| \times |D_n^j|}.$$

Applying conditions (C2), (2.1) and Kronecker's lemma, we conclude that

$$|D_n| \mathrm{cov}\{\widehat{C}_n(\mathbf{h}_i, u_i), \widehat{C}_n(\mathbf{h}_j, u_j)\} \to \sum_{\mathbf{s} \in \mathbb{Z}^2} \sum_{t \in \mathbb{Z}} \mathrm{cov}\{Z(\mathbf{0},0) Z(\mathbf{h}_i, u_i), Z(\mathbf{s},t) Z(\mathbf{s} + \mathbf{h}_j, t + u_j)\}.$$

Let $\sigma^2 = \sum_{\mathbf{s} \in \mathbb{Z}^2} \sum_{t \in \mathbb{Z}} \operatorname{cov}\{Z(\mathbf{0},0)Z(\mathbf{h}_1, u_1), Z(\mathbf{s},t)Z(\mathbf{s} + \mathbf{h}_2, t + u_2)\}$ and $A_n = \sqrt{|D_n|}\{\widehat{C}_n(\mathbf{h}, u) - C_n(\mathbf{h}, u)\}$. By Lemma A.2, $A_n \xrightarrow{d} \mathrm{N}(0, \sigma^2)$. We apply the Cramér–Wold device to prove the joint normality. □

**Proof of Theorem 2.** See the proof in Li *et al.* (2007). □

**Proof of Theorem 3.** This is a special case of Lemma A.4, setting $p = 2$ and $q = 1$. □

**Proof of Theorem 4.** Asymptotic normality of $\widehat{\mathbf{G}}_n$ follows immediately from Lemma A.4, setting $p = 3$, $q = 0$. □

## Acknowledgements

The National Center for Atmospheric Research is sponsored by the National Science Foundation. Genton acknowledges partial support from National Science Foundation grants DMS-0504896 and CMG ATM-0620624. The authors thank the Editor, the Associate Editor and a referee for constructive suggestions that have improved the content and presentation of this article.